\documentclass[12pt,a4paper]{article}
\usepackage{arxiv}
\usepackage{graphicx}
\usepackage[T1]{fontenc}
\usepackage[utf8]{inputenc}
\usepackage[english]{babel}
\usepackage{amsmath,amssymb,amsfonts,amsthm}
\usepackage{hyperref}
\usepackage{cite}
\usepackage{xcolor} % revision highlighting

\AtBeginDocument{\setlength{\headheight}{15pt}\raggedbottom}
\newtheorem*{remark}{Remark}

\title{Fast Eighth-Order Pad\'e Schemes Based on Chebyshev Polynomials 
for the Direct Zakharov--Shabat Problem}
\author{Sergey Medvedev$^{1,2,*}$, Igor Chekhovskoy$^{2,1}$, Irina Vaseva$^{1,2}$, Mikhail Fedoruk$^{2,1}$\\
$^{1}$ Federal Research Center for Information and Computational Technologies,\\ 
Novosibirsk 630090, Russia,\\
$^{2}$ Novosibirsk State University, Novosibirsk 630090, Russia,\\
* Corresponding author: medvedev@ict.nsc.ru}
\date{}

\begin{document}
	
\maketitle
	
\begin{abstract}
In this work, we construct fast eighth-order Pad\'e schemes for the direct Zakharov--Shabat scattering problem. The schemes are based on an eighth-order exponential integrator obtained from the Magnus expansion. A direct extension of the conventional fast Pad\'e representation to the eighth-order case leads to insufficiently accurate fast variants in the considered tests. To overcome this difficulty, we reformulate the spectral dependence on a finite real interval using the Joukowski mapping and represent the local numerators and denominators in the Chebyshev polynomial basis. This makes it possible to construct the global transition matrix by a fast product tree while retaining a compact polynomial representation. Numerical experiments for chirped hyperbolic secant potentials with both signs of dispersion show that the proposed Chebyshev-based fast schemes substantially improve the accuracy of the corresponding direct fast variants and can be used for efficient computation of the continuous nonlinear spectrum.
\end{abstract}

	\section{Introduction}
	
	The direct Zakharov--Shabat (ZS) problem is the first stage of the inverse scattering transform method, which is also known as the nonlinear Fourier transform (NFT). The NFT is a technique for solving the Cauchy problem for nonlinear evolution equations. It is based on the connection between a nonlinear equation and the scattering data of a family of auxiliary linear differential operators, which makes it possible to reconstruct the evolution of the solution of the nonlinear equation from the evolution of the scattering data. The scattering data are also called spectral data or the nonlinear spectrum of the problem. At the first stage of the NFT, for a known potential $q(z_{0}, t)$, the direct scattering problem is solved and the spectral data are found. Next, the evolution of the spectral data is determined using elementary formulas. At the third stage, the inverse problem is solved, i.e., the potential $q(z, t)$ at the point $z$ is determined from the known scattering data.
	
	The NFT was first proposed in 1967 by Gardner, Greene, Kruskal, and Miura for the Korteweg--de Vries equation~\cite{Gardner1967}. In 1971, Zakharov and Shabat~\cite{ZakharovShabat1972}, using the approach first introduced by Lax~\cite{Lax}, showed that the NFT is applicable to solve the nonlinear Schr\"odinger equation (NLSE)
	\begin{equation} \label{nlse}
		i\frac{\partial q}{\partial z} +\frac{\sigma}{2} \frac{\partial ^{2} q}{\partial t^{2} } +|q|^{2} q=0,
	\end{equation}
	which describes the propagation of pulses in an optical fiber~\cite{Hasegawa1973a,Hasegawa1973b,Agrawal}.
	In Eq.~\eqref{nlse}, $q(t,z)$ is a %\textcolor{red}
    {slowly varying} complex optical field envelope, $z$ is the distance along the optical fiber, and $t$ is a time variable. The parameter $\sigma$ denotes anomalous ($\sigma=1$) or normal dispersion ($\sigma=-1$) in the fiber.
	
	Using the NFT, the complex nonlinear dynamics of the NLSE solution can be reduced to a simple evolution in the nonlinear ZS spectrum. %\textcolor{red}
    {In telecommunications applications,} the NFT is more often used as a signal processing and encoding tool~\cite{Yousefi2014I}.
	At present, the main goal of the numerical solution to the direct ZS problem is to design fast algorithms with low computational complexity $\mathcal{O}(N_t\log^2{N_t})$, where $N_t$ is the number of signal samples~\cite{Wahls2013,Wahls2015a,Wahls2018}. The numerical NFT literature includes layer-peeling and spectral methods~\cite{Yousefi2014II}, the Boffetta--Osborne discretization~\cite{Boffetta1992a}, Ablowitz--Ladik-type discretizations~\cite{Ablowitz1975NonlinearEquations,Ablowitz1976AScattering}, Volterra and integral-equation approaches for scattering-data computation~\cite{Fermo2015ZS,Frumin2015_TIB}, extensions of scattering algorithms to vector Manakov systems~\cite{Frumin2020}, contour-integral and hybrid methods for eigenvalue detection~\cite{Vasylchenkova2018a,Vasylchenkova2019a,Aref2019}, numerical direct scattering for large wave packets and breather-type fields~\cite{Mullyadzhanov2019,Mullyadzhanov2023Breathers}, spectral-parameter power-series representations for direct and inverse ZS scattering~\cite{Kravchenko2025SPPS}, as well as studies of numerical errors in the direct scattering transform~\cite{Gelash2020_PRE}. The idea of fast methods for the direct Zakharov--Shabat problem is to represent the local transition matrix of the system as a rational or polynomial function of the spectral parameter. Then fast polynomial arithmetic algorithms can be used to compute the product of local transition matrices. If the transition matrix of the ZS system has an exponential form, it can be transformed into a suitable rational form by using Pad\'e approximation or exponential splitting schemes. Such methods are combined under the general name of fast nonlinear Fourier transform (FNFT).
	
	Fast schemes based on commutator-free quasi-Magnus (CFQM) exponential integrators were presented in~%\textcolor{red}
    {\cite{Chimmalgi2019}}. In addition to fourth-order schemes, the authors proposed a sixth-order CFQM scheme with complex coefficients, as well as a sixth-order scheme obtained from a fourth-order CFQM scheme using Richardson extrapolation. A distinctive feature of CFQM-based schemes is the need to use an irregular grid, while in telecommunication applications the signal is usually specified in a tabular form on an equidistant grid. Higher-order fast algorithms were also constructed using multistep discretizations~\cite{Vaibhav2018}, Runge--Kutta methods~\cite{vaibhav2019efficient}, and Chebyshev-polynomial representations for the continuous and discrete spectral data~\cite{Vaibhav2019Chebyshev,Vaibhav2019ChebyshevSpectral,Vaibhav2019ChebyshevNorming}. Fast sixth-order schemes for solving the direct Zakharov--Shabat problem based on the Magnus expansion~\cite{Magnus1954,Blanes2009}, the generalized Cayley transform, and diagonal Pad\'e approximation were proposed in~\cite{MedvedevPade}. In~\cite{Medvedev2020_OL,Medvedev2024splitting} exponential splitting schemes were examined to construct FNFT schemes based on the three-exponential fourth-order scheme TES4.
	
	Among the conventional high-order schemes, one can note the works~\cite{blanes2006fourth,Blanes2017}, which present sixth-order schemes based on CFQM exponential integrators. In~\cite{alvermann2011high}, eighth-order CFQM schemes for solving a general non-autonomous system are presented. It should also be noted that for a more balanced solution of the direct ZS problem, it is important to preserve the quadratic invariant for the real spectral parameter. Conservation of quadratic invariants is important for long-time numerical integration~\cite{Hairer1987}. Fast conservative schemes include the fourth-order CFQM scheme $\text{CF}^{[4]}_2$~\cite{Chimmalgi2018}, a scheme based on TES4 and the Suzuki decomposition~\cite{Medvedev2020_OL}, and sixth-order schemes based on the generalized Cayley transform and the diagonal Pad\'e approximation~\cite{MedvedevPade}.
	
	An approach to constructing high-order schemes for implementing the inverse NFT was recently proposed in~\cite{MedvedevHGTIB}. Using Gregory's formula with an appropriate number of coefficients, it is possible to construct a scheme with up to the ninth order of accuracy~\cite{Medvedev2025inverse9}. This prompted an attempt to also improve the order of accuracy for solving the direct ZS problem for a consistent implementation of the direct and inverse NFT. In the present work, we consider an eighth-order exponential scheme, denoted below as ES8, based on the Magnus expansion~\cite{Magnus1954}. The diagonal Pad\'e approximation is employed to convert the local transition matrix into rational form.
	
	A direct fast implementation of the resulting eighth-order Pad\'e schemes can be constructed by following the same general idea as in lower-order methods. However, numerical experiments show that the corresponding direct fast variants are not sufficiently accurate for the considered continuous-spectrum tests. Therefore, we propose an alternative representation based on a finite spectral interval, the Joukowski mapping, and Chebyshev polynomials. %\textcolor{red}
    {The Joukowski mapping connects the finite interval with the unit circle, while the Chebyshev basis provides a compact representation for the fast product-tree multiplication.} Chebyshev polynomials were previously applied to FNFT in~\cite{Vaibhav2019Chebyshev,Vaibhav2019ChebyshevSpectral,Vaibhav2019ChebyshevNorming}. In contrast to those works, the local transition matrix considered here is generated by the diagonal Pad\'e approximation of an eighth-order exponential scheme obtained from the Magnus expansion.
	
	The paper is organized as follows. Section~2 recalls the direct Zakharov--Shabat problem and the notation used for the scattering coefficients. Section~3 describes the construction of the fast eighth-order Pad\'e schemes in the Chebyshev polynomial basis. Section~4 presents numerical experiments for chirped hyperbolic secant potentials with both signs of dispersion. The Appendix contains the explicit formulas required for the implementation of the eighth-order exponential scheme.
	
	\section{Zakharov--Shabat Problem}
	The direct ZS problem for the NLSE~\eqref{nlse} with the complex spectral parameter $%\textcolor{red}
    {\zeta=\xi+i\eta}$ can be written as an evolutionary system
	\begin{equation}\label{psit}
		\frac{d{\Psi}(t)}{dt}=Q(t){\Psi}(t),\quad 
		Q(t)=\left(\begin{array}{cc}-i\zeta&q\\r&i\zeta\end{array}\right),\quad \Psi(t) = \left(\begin{array}{c} \psi_1\\ \psi_2\end{array}\right),
	\end{equation}
	where $q=q(t,z_0)$ is the initial field for the NLSE at the point $z=z_0$, $r = -\sigma q^*$;
	$\Psi(t)$ is a wave function. 
	Here, $z_0$ plays the role of a parameter and will be omitted. It is assumed that $q(t)$ decays rapidly when $t\to \pm \infty$.
	
	The main goal of solving the direct ZS problem~\eqref{psit} is to find the nonlinear spectrum for a given signal $q(t)$. In general, the 
	nonlinear 
	spectrum consists of a continuous and a discrete spectrum and is described by the scattering coefficients $a(\zeta)$ and $b(\zeta)$. They are 
	defined 
	for the real spectral parameter  $\zeta=\xi$ as follows
	\begin{equation}\label{ab}
		a(\xi)=\lim_{t\to\infty}\,\psi_1(t,\xi)\,e^{i\xi t},\quad b(\xi)=\lim_{t\to\infty}\,\psi_2(t,\xi)\,e^{-i\xi t},
	\end{equation}
	where the wave function $\Psi(t)$ is determined at $t\to-\infty$ by the Jost functions
	\begin{equation}\label{psi0}
		\Psi =\left(
		\begin{array}{c}
			\psi_{1}\\\psi_{2}
		\end{array}
		\right) = \left(
		\begin{array}{c}
			e^{-i\zeta t}\\0
		\end{array}
		\right)[1+o(1)],\quad t\to-\infty.
	\end{equation}
	Similarly, one can write the Jost functions for $t\to +\infty$ and define corresponding scattering coefficients at $t\to -\infty$.
	
	The scattering coefficients~\eqref{ab} determine the reflection coefficient
	$%\textcolor{red}
    {\rho(\xi)}={b(\xi)}/{a(\xi)}$, $\xi\in\mathbb{R}$, which describes the continuous spectrum of the ZS problem.
	
	For anomalous dispersion, the discrete spectrum can exist. To define the discrete spectrum, the scattering coefficients are continued into the 
	upper 
	complex half-plane. Complex eigenvalues of the ZS problem are found as zeros of $a(\zeta)$, where $\zeta$ is a complex number with a positive 
	imaginary 
	part. The discrete spectrum consists of the eigenvalues $\zeta_k$ and associated phase coefficients
	$
	r_k=\left.b(\zeta)/a'(\zeta)\right|_{\zeta=\zeta_k},\ \mbox{where}\  a'(\zeta)=da(\zeta)/d\zeta.
	$
	
	The ZS system~\eqref{psit} conserves the quadratic invariant $|\psi_1(\xi)|^2+\sigma|\psi_2(\xi)|^2 \equiv 1$ for real values of the spectral 
	parameter, %\textcolor{red}
    {which implies}
	\begin{equation}\label{H}
		H(\xi) = |a(\xi)|^2+\sigma|b(\xi)|^2 \equiv 1.
	\end{equation}
	
	%he main features of the computational problem can be found in~\cite{Medvedev2020_OE}. Briefly, the unknown function $\Psi$ must be calculated on 
	%a 
	%uniform grid; Dahlquist's Second Barrier restricts the application of multistep methods ~\cite{Dahlquist1963, Hairer1987}; matrix exponentials 
	%can 
	%be 
	%easily calculated for matrices 2x2; it is necessary to solve the ZS system for a large number of values of spectral parameters $\zeta$ at a 
	%fixed 
	%potential~$q(t)$.
	
	The numerical implementation of the continuous function~$q(t)$ is a discrete function~$q_n=q(t_n)$, which is defined at the nodes $t_n$ of the uniform grid with the step size $\tau$. Since we are considering a finite time interval, we solve the problem on $[-L,L]$ with the total number of time steps equal to $N_t$, the grid step size $\tau=2L/N_t$, and $t_n=-L+\tau n$, where $n=0,\ldots,N_t$. The original 
	system~\eqref{psit} on the interval $(t_n-\tau/2, t_n+\tau/2)$ is replaced with the approximate system with constant coefficients
	\begin{equation}\label{Tn}
		\Psi(t_n+\tau/2)=T_n\Psi(t_n-\tau/2),
	\end{equation}
	where $T_n$ is a local transition matrix.
	
	The transition matrix $T_n$ can be defined using the Magnus expansion \cite{Magnus1954, Blanes2009}, which gives the asymptotic 
	representation 
	of the fundamental solution of the system~\eqref{psit}
	$$\Psi(t)=U(t,0)\Psi(0),\quad U(t,0)=e^{\Omega(t)},\quad \Omega(t)=%\textcolor{red}
    {\sum\limits_{k=1}^\infty \Omega_k(t)}.$$
	First terms of the Magnus expansion have the form
	$$\Omega_1(t)=\int\limits_{0}^{t}dt_1 Q(t_1),\quad
	\Omega_2(t)=\frac{1}{2}\int\limits_{0}^{t}dt_1\int\limits_{0}^{t_1}dt_2\,[Q(t_1), Q(t_2)],$$
	$$\Omega_3(t)=\frac{1}{6}\int\limits_{0}^{t}dt_1\int\limits_{0}^{t_1}dt_2\int\limits_{0}^{t_2}dt_3\,\left([Q(t_1),[Q(t_2), 
	Q(t_3)]]+[Q(t_3),[Q(t_2),Q(t_1)]]\right).$$
	Square brackets $[A,B]$ are the matrix commutator of $A$ and $B$.
	
	If the matrix $Q(t)$ can be represented as a Taylor series with respect to a small parameter $\tau$
	\begin{equation}
		Q(t+\tau)=\sum\limits_{k=0}^\infty\,\frac{\tau^k}{k!}Q^{(k)}(t),\quad Q^{(k)}(t)=\frac{d^kQ(t)}{dt^k},
	\end{equation}
	then, substitution of this series into the Magnus formula with integration from $t-\tau/2$ to $t+\tau/2$ gives an approximation of the 
	fundamental 
	solution of the problem with the required order of accuracy in the small parameter $\tau$.
	Here, we consider the approximation up to the eighth order in $\tau$. Then, the transition matrix of the ZS system has the form
	\begin{equation}\label{ES8}
		T^{(ES8)}=e^{Z(t)},\quad Z(t)=\tau \hat{Z}_1(t)+\tau^3\hat{Z}_3(t)+\tau^5\hat{Z}_5(t) +\tau^7 \hat{Z}_7(t).
	\end{equation}
	Expressions for the matrices $\hat{Z}_1$, $\hat{Z}_3$, $\hat{Z}_5$, $\hat{Z}_7$ for the system~\eqref{psit} are determined by the formulas~\eqref{Z1eZ3e},~\eqref{Z5e}, and~\eqref{Z7} presented in the Appendix. Some additional computational details can also be found in the Appendix.

	\section{Fast scheme implementation}

	When constructing a numerical scheme for the ZS problem, a local transition matrix~\eqref{Tn}
	is formed at each step of the $t$-discretization. To compute the scattering coefficients  $a(\zeta)$ and $b(\zeta)$ one needs to construct the 
	global 	transition 	matrix
	\begin{equation}
		T(\zeta) = T_{N_t}(\zeta)T_{N_t-1}(\zeta)\ldots T_1(\zeta).
	\end{equation}
	If the matrices $T_n(\zeta)$ have an arbitrary dependence on $\zeta$, then computing  
	$T(\zeta)$ for a large number of spectral parameters requires a large number of independent  
	solutions of the initial problem. Therefore, to construct a fast algorithm, one needs to ensure that $T_n(\zeta)$  
	depends polynomially or rationally on $\zeta$. This is precisely what Pad\'e approximations provide.
	This makes it possible to represent the local transition matrix as a rational function of the spectral parameter. Consequently, the product of 
	the local transition matrices over the entire grid is reduced to multiplying matrix polynomials $N_n(\zeta)$ and scalar denominators 
	$d_n(\zeta)$: 
	\begin{equation}\label{LocalT}
		T_n(\zeta)
		=
		\frac{1}{d_n(\zeta)}N_n(\zeta),
	\end{equation}
	Hence, fast polynomial arithmetic 
	algorithms become applicable for computing the global transition matrix
	\begin{equation}\label{globalT}
		T(\zeta)
		=
		\frac{1}{D(\zeta)}N(\zeta),
	\end{equation}
	where
	\begin{equation}\label{product}
		N(\zeta)
		=
		N_{N_t}(\zeta)N_{N_t-1}(\zeta)\ldots N_1(\zeta),\quad
		D(\zeta)
		=
		d_{N_t}(\zeta)d_{N_t-1}(\zeta)\ldots d_1(\zeta).
	\end{equation}
	The numerator and denominator can be constructed independently.
	
	For example, for the diagonal Pad\'e approximation of degree $s$ $E_s(z)	=	P_s(z)/P_s(-z)$
	the local transition matrix has the form
	\begin{equation}
		T_n(\zeta)
		=
		E_s(Z_n(\zeta))
		=
		P_s(Z_n(\zeta))P_s(-Z_n(\zeta))^{-1}.
	\end{equation}
	For the ZS system, this expression can be reduced to a rational form
	\begin{equation}\label{RationalForm}
		T_n(\zeta)
		=
		f_n(\zeta)E+g_n(\zeta)Z_n(\zeta),
	\end{equation}
	where $E$ is the identity matrix, $Z_n(\zeta)$ is the local matrix of the exponential scheme, and
	$f_n$ and $g_n$ are rational functions of the spectral parameter.
	
	In conventional fast schemes, the Cayley transform is usually used, which maps the real axis of the spectral parameter to the unit circle. After such a transformation, fast polynomial algorithms can be employed. This approach is efficient for a number of lower-order schemes, including the sixth-order Pad\'e schemes based on the generalized Cayley transform~\cite{MedvedevPade}. For the eighth-order Pad\'e schemes considered in this work, however, the direct fast variants obtained in this way were found to be insufficiently accurate in the numerical tests presented below. This motivates a different representation of the spectral dependence.
	
	In practice, the continuous spectrum is typically needed not on the entire real axis,
	but on a specified finite interval $\zeta\in[\zeta_L,\zeta_R].$
	Hence, it is natural to abandon the mapping of the entire real axis and switch to a finite
	interval:
	\begin{equation}
		\zeta=c+Hx,\quad	c=\frac{\zeta_R+\zeta_L}{2},\quad
		H=\frac{\zeta_R-\zeta_L}{2}%\textcolor{red}
        {,\quad} x\in[-1,1].
	\end{equation}
	In the symmetric case $c=0$, we have $\zeta=Hx$.
	
	Thus, the spectral problem is rewritten in the variable $x$, which lies on a finite
	interval. This removes the problem of compressing infinity, but now it is
	necessary to choose a suitable basis 	for polynomials in $x$.
	
	The finite interval $[-1,1]$ is naturally connected to the unit circle via the Joukowski
	transform:
	\begin{equation}\label{Zhuk}
		x=\frac{1}{2}\left(w+w^{-1}\right).
	\end{equation}
	If $w=e^{i\theta}$, then $x=\cos\theta$. 
	Consequently,
	\begin{equation}
		\zeta=c+H\cos\theta.
	\end{equation}
	The circle is thereby mapped onto the finite interval $[\zeta_L,\zeta_R]$. This is the main
	difference from the Cayley transform: we do not attempt to cover the entire real axis with a single
	rational mapping, but instead work directly on the finite spectral interval.
	
	If a local polynomial in $x$ has degree $K$,
	\begin{equation}
		p(x)=\sum_{k=0}^{K}p_kx^k,
	\end{equation}
	then substituting Eq.~\eqref{Zhuk} we obtain a Laurent polynomial in $w$ with support from $-K$ to $K$:
	\begin{equation}
		p\left(\frac{w+w^{-1}}{2}\right)
		=
		\sum_{k=-K}^{K}\widehat{p}_k w^k.
	\end{equation}
	If one then clears the negative powers by multiplying by $w^K$, an ordinary polynomial
	of degree $2K$ is obtained.
	
	\begin{remark}
		This is precisely where the danger of an inefficient implementation arises. If the Joukowski
		transform is used directly and the cleared Laurent polynomial is stored, the local degree effectively
		doubles. For example, an object of degree $30$ in $x$ becomes an ordinary polynomial of degree $60$
		in $w$. This increases the constant in the fast algorithm.
	\end{remark}
	
	Consequently, the Joukowski transform itself is useful as a geometric idea, but computationally
	it is better not to switch to cleared Laurent polynomials. The proper basis for the variable
	$x=\cos\theta$ is the Chebyshev polynomial basis.

	The Chebyshev polynomials of the first kind are defined by the formula
	\begin{equation}
		T_k(x)=\cos(k\arccos x).
	\end{equation}
	For $x=\cos\theta$, we have $T_k(\cos\theta)=\cos(k\theta)$. 
	If $w=e^{i\theta}$, then
	\begin{equation}
		\cos(k\theta)
		=
		\frac{1}{2}\left(e^{ik\theta}+e^{-ik\theta}\right)
		=
		\frac{1}{2}\left(w^k+w^{-k}\right).
	\end{equation}
	From this, the fundamental relation follows:
	\begin{equation}\label{ChebRelation}
		T_k(x)
		=
		\frac{1}{2}\left(w^k+w^{-k}\right),
		\quad
		x=\frac{1}{2}\left(w+w^{-1}\right).
	\end{equation}
	
	The values of the Chebyshev polynomials can be computed recursively:
	\begin{equation}
		T_0(x)=1,
		\quad
		T_1(x)=x,\quad
		T_{k+1}(x)=2xT_k(x)-T_{k-1}(x).
	\end{equation}
	If values are needed simultaneously at a large number of spectral nodes, one can use the connection
	with the relation~\eqref{ChebRelation}.

	Thus, Chebyshev polynomials are a natural way to work with the Joukowski variable without
	doubling the number of stored coefficients. If
	\begin{equation}
		p(x)=\sum_{k=0}^{K}%\textcolor{red}
        {c_k}T_k(x),
	\end{equation}
	then storing such a polynomial requires $K+1$ coefficients, rather than the $2K+1$ coefficients
	of the Laurent polynomial.

	%\textcolor{red}
    {In the implementation, the coefficients $c_k$ are obtained directly from the power-basis coefficients before the product tree is built. If $p(x)=\sum_{n=0}^{K}a_nx^n$, then}
	\begin{equation}\label{PowerToCheb}
		%\textcolor{red}
        {c_k=\sum_{\substack{n=k\\ n-k\ {\rm even}}}^{K}\alpha_{n,k}a_n,\qquad
		\alpha_{n,0}=2^{-n}\binom{n}{n/2},\qquad
		\alpha_{n,k}=2^{1-n}\binom{n}{(n-k)/2}\ \ (k>0).}
	\end{equation}
	%\textcolor{red}
    {Thus, the change of basis is exact and uses only precomputed binomial factors; no interpolation or discrete cosine transform is required.}
	
	\begin{remark}
		Geometrically, the Joukowski transform explains why the variable $x=\cos\theta$ arises.
		Algebraically, Chebyshev polynomials provide a compact basis for this variable.
		Therefore, it is natural to call the final scheme a scheme in the Chebyshev polynomial basis,
		rather than simply a Joukowski scheme.
	\end{remark}
	
	After switching from $\zeta$ to $x$, the local transition matrix is written in the form~\eqref{LocalT}.
	The numerator and denominator are expanded in the Chebyshev basis:
	\begin{equation}
		N_n(x)
		=
		\sum_{k=0}^{K}N_{n,k}T_k(x),\quad
		d_n(x)
		=
		\sum_{k=0}^{K_d}d_{n,k}T_k(x).
	\end{equation}
	Here $N_{n,k}$ are $2\times 2$ matrices, and $d_{n,k}$ are scalar coefficients.
	
	The degrees $K$ and $K_d$ are determined only by the local scheme. If the %\textcolor{red}
    {local matrix}
	$Z_n(\zeta)$ has a polynomial dependence on $\zeta$ of degree $p$, then after the
	Pad\'e approximation of degree $s$, the %\textcolor{red}
    {maximum degree of the numerator and denominator in $\zeta$ is $K=2ps$}.
	The factor $2$ arises from the structure of the rational function of a matrix and the transition to
	functions of the quadratic combination of the matrix components.
	
	For eighth-order schemes, the degree of the local matrix in the spectral parameter is
	$p_{\mathrm{ES8}}=5$, and therefore $K_{\mathrm{ES8},s}=10s$.
	For example, for the ES8 scheme with Pad\'e approximation of degree $3$, we have
	$K_{\mathrm{ES8},3}=30$.

	To construct the global transition matrix~\eqref{globalT}, it is necessary to multiply polynomials in the Chebyshev basis.
	The fundamental formula is
	\begin{equation}
		T_j(x)T_k(x)
		=
		\frac{1}{2}
		\left(
		T_{j+k}(x)+T_{|j-k|}(x)
		\right).
	\end{equation}
	Let
	\begin{equation}
		p(x)=\sum_{j=0}^{m}p_jT_j(x),
		\quad
		q(x)=\sum_{k=0}^{n}q_kT_k(x),
	\end{equation}
	then
	\begin{equation}
		r(x)=p(x)q(x)
		=
		\sum_{\ell=0}^{m+n}r_\ell T_\ell(x).
	\end{equation}
	The coefficients $r_\ell$ are obtained from two types of contributions:
	$j+k=\ell$ and
	$|j-k|=\ell$.
	The first contribution has the form of an ordinary convolution, the second one is a correlation term. Therefore,
	multiplication in the Chebyshev basis can be implemented via fast convolutions. In a practical
	implementation, it is important to reuse temporary arrays and accumulate the result directly
	into the appropriate element of the matrix polynomial.

	%\textcolor{red}
    {In the present implementation, one scalar Chebyshev product is evaluated by two FFT-based ordinary polynomial products corresponding to the convolution and correlation contributions in the formula above. The correlation term is obtained by reversing one coefficient array. For polynomials of degree $L$, each ordinary polynomial product uses two forward FFTs and one inverse FFT; hence one scalar Chebyshev product uses four forward FFTs and two inverse FFTs, followed by $O(L)$ coefficient recombination. Its asymptotic complexity remains $O(L\log L)$, although with a larger constant than ordinary polynomial multiplication. The standard $2\times2$ matrix-polynomial product requires eight such scalar Chebyshev products.}
	
	For matrix polynomials
	\begin{equation}
		A(x)=
		\begin{pmatrix}
			A_{11}(x) & A_{12}(x)\\
			A_{21}(x) & A_{22}(x)
		\end{pmatrix},
		\quad
		B(x)=
		\begin{pmatrix}
			B_{11}(x) & B_{12}(x)\\
			B_{21}(x) & B_{22}(x)
		\end{pmatrix}
	\end{equation}
	the product has the standard form
	$C(x)=A(x)B(x)$, where $C_{ij} = A_{i1}B_{1j}+A_{i2}B_{2j}$.
	Each product on the right-hand side is a product of Chebyshev polynomials.
	
	\begin{remark}
		In the implementation, one should not create a separate temporary polynomial for each product
		$A_{ik}B_{kj}$. It is better to accumulate the contribution directly into the corresponding element
		$C_{ij}$. This does not change the mathematical scheme, but reduces the number of memory allocations
		and copies.
	\end{remark}

	To construct the global transition matrix~\eqref{globalT}, both products in~\eqref{product} are built using a product tree. At the first level, 
	adjacent 
	factors are multiplied:
	\begin{equation}
		N_{2j,2j-1}(x)=N_{2j}(x)N_{2j-1}(x).
	\end{equation}
	At the next level, the already obtained products are multiplied. This continues until
	only one global matrix polynomial remains.
	
	If at some level of the tree the degree of the polynomials is of order $L$, then the fast
	polynomial multiplication has complexity
	$O(L\log L)$.
	At this level, the number of polynomials is inversely proportional to $L$, so the total
	cost of the level is
	$O(N_t\log N_t)$.
	Since the number of levels is $O(\log N_t)$, the total complexity of building the product tree is
	$O(N_t\log^2N_t)$.

	%\textcolor{red}
    {At every level of the tree, multiplying Chebyshev polynomials of degrees $m$ and $n$ produces a Chebyshev polynomial of degree $m+n$, stored using only $m+n+1$ coefficients. The conversion to a Laurent representation is performed only after the complete product tree has been constructed. If the final Chebyshev polynomial is $p(x)=\sum_{k=0}^{K_{\rm g}}c_kT_k(x)$, relation~\eqref{ChebRelation} gives}
	\begin{equation}\label{ChebToLaurent}
		%\textcolor{red}
        {\widehat p_0=c_0,\qquad
		\widehat p_{-k}=\widehat p_k=\frac{1}{2}c_k,\quad 1\le k\le K_{\rm g}.}
	\end{equation}
	%\textcolor{red}
    {Hence the final Laurent polynomial contains the powers $w^{-K_{\rm g}},\ldots,w^{K_{\rm g}}$ and has $2K_{\rm g}+1$ coefficients. After multiplication by $w^{K_{\rm g}}$, it is equivalent to an ordinary polynomial of degree $2K_{\rm g}$. The factor-of-two increase therefore appears only at this final conversion, not during the product-tree multiplication. For example, the local ES8 Pad\'e scheme with $s=3$ has working Chebyshev degree $K=30$ and stores $31$ coefficients, whereas its direct Laurent/Joukowski representation has $61$ coefficients and is equivalent to an ordinary polynomial of degree $60$. The conventional Cayley-based ES8 Pad\'e scheme of the same Pad\'e degree has local polynomial degree $30$.}
	
	\begin{remark}
		Here $N_t$ is the number of signal samples. The degree of the local polynomial is considered a fixed
		constant that depends only on the chosen numerical scheme and the degree of the Pad\'e approximation.
	\end{remark}
	
	After building the product tree, the global transition matrix has the form
	\begin{equation}
		T(x)
		=
		\frac{1}{D(x)}
		\begin{pmatrix}
			N_{11}(x) & N_{12}(x)\\
			N_{21}(x) & N_{22}(x)
		\end{pmatrix}.
	\end{equation}
	If the first column of the transition matrix is used, the spectral coefficients are determined by
	the formulas
	\begin{equation}
		a(\zeta)
		=
		\frac{N_{11}(x)}{D(x)},
		\quad
		b(\zeta)
		=
		\frac{N_{21}(x)}{D(x)},
		\quad
		x=\frac{\zeta-c}{H}.
	\end{equation}
	Thus, the new approach does not change the definition of the coefficients $a(\zeta)$ and $b(\zeta)$.
	Only the method of constructing the global polynomials in the spectral parameter changes.

	\section{Numerical Results}
	
In this section, we compare a non-fast eighth-order Pad\'e scheme used as a reference, direct fast Pad\'e variants, and the proposed Chebyshev-based fast schemes.
The goal of the numerical experiments is to verify whether the Chebyshev representation makes it possible to obtain practical fast variants of the eighth-order scheme for computing the continuous nonlinear spectrum.
	
	As a model potential, we use the chirped hyperbolic secant profile
	\begin{equation}
		q(t)=A[\operatorname{sech}(t)]^{1+iC},
	\end{equation}
	where $A=5.2$ and $C=4$. This test is more demanding than the unchirped case and was also used in our previous numerical comparisons of NFT schemes~\cite{Medvedev2020_OE,Medvedev2024splitting}. We consider both signs of dispersion, $\sigma=1$ and $\sigma=-1$. The continuous spectrum is computed on the finite interval $\xi\in[-20,20]$. In the main experiments, the number of spectral points is fixed as $N_\xi=2^{15}$, while the number of time steps $N_t$ varies.
	
	The numerical error is measured by the normalized mean squared error
	\begin{equation}\label{NMSE}
		\operatorname{NMSE}[\phi]
		=
		\frac{1}{N_\xi}
		\sum_{j=1}^{N_\xi}
		\frac{|\phi^{\mathrm{comp}}(\xi_j)-\phi^{\mathrm{exact}}(\xi_j)|^2}{|\phi_0(\xi_j)|^2},
	\end{equation}
	where $\phi$ denotes either $a(\xi)$ or $b(\xi)$, and
	\begin{equation}
		\phi_0(\xi_j)=
		\begin{cases}
			\phi^{\mathrm{exact}}(\xi_j), & |\phi^{\mathrm{exact}}(\xi_j)|>1,\\
			1, & |\phi^{\mathrm{exact}}(\xi_j)|\le 1.
		\end{cases}
	\end{equation}
	For each scheme, we plot the error as a function of the time grid size and as a function of the execution time. This presentation makes it possible to compare not only the asymptotic accuracy of the schemes, but also their practical efficiency.
	
	\begin{figure}[t]
		\centering
		\includegraphics[width=0.8\textwidth]{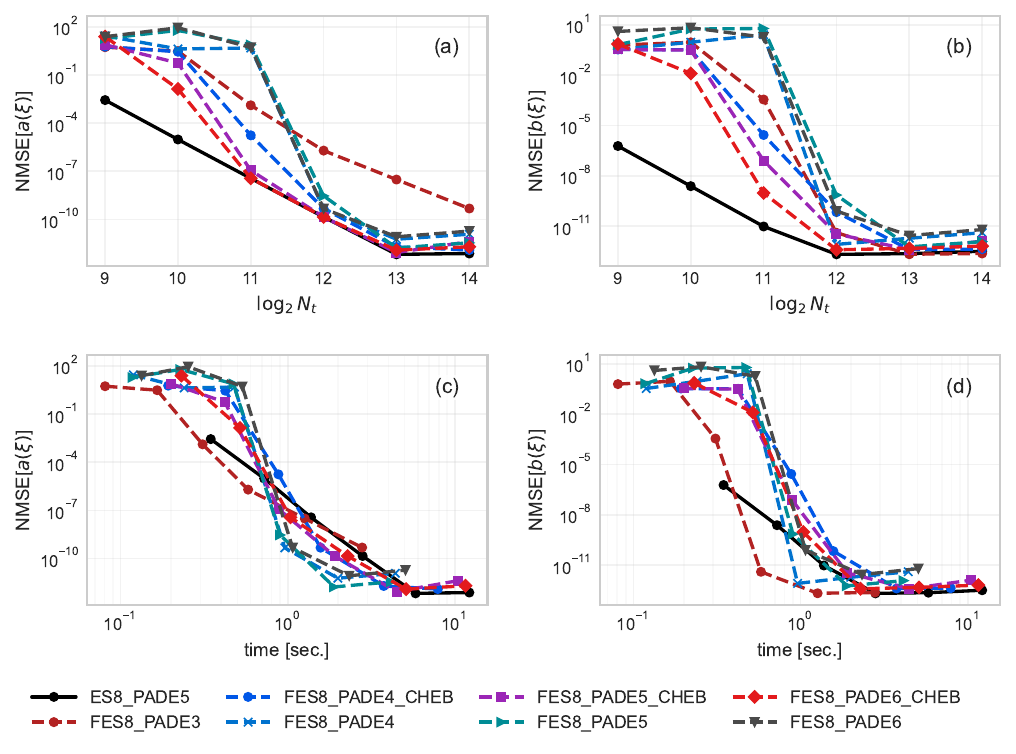}
		\caption{Continuous-spectrum errors for the chirped hyperbolic secant potential with $A=5.2$, $C=4$, and anomalous dispersion $\sigma=1$. The upper panels show $\operatorname{NMSE}[a(\xi)]$ and $\operatorname{NMSE}[b(\xi)]$ as functions of the time grid size. The lower panels show the same errors as functions of the execution time. The spectral interval is $\xi\in[-20,20]$ and $N_\xi=2^{15}$.}
		\label{fig:es8_sigma_plus}
	\end{figure}

%	Figure~\ref{fig:es8_sigma_plus} shows the results for anomalous dispersion, $\sigma=1$. The direct fast variants reproduce the general convergence trend, but their accuracy is limited compared with the corresponding non-fast Pad\'e schemes. The Chebyshev-based fast variants substantially reduce this gap. In particular, for sufficiently resolved time grids they provide noticeably smaller errors in both $a(\xi)$ and $b(\xi)$ than the direct fast variants, while preserving the fast product-tree structure of the algorithm.
	
	\begin{figure}[t]
		\centering
		\includegraphics[width=0.8\textwidth]{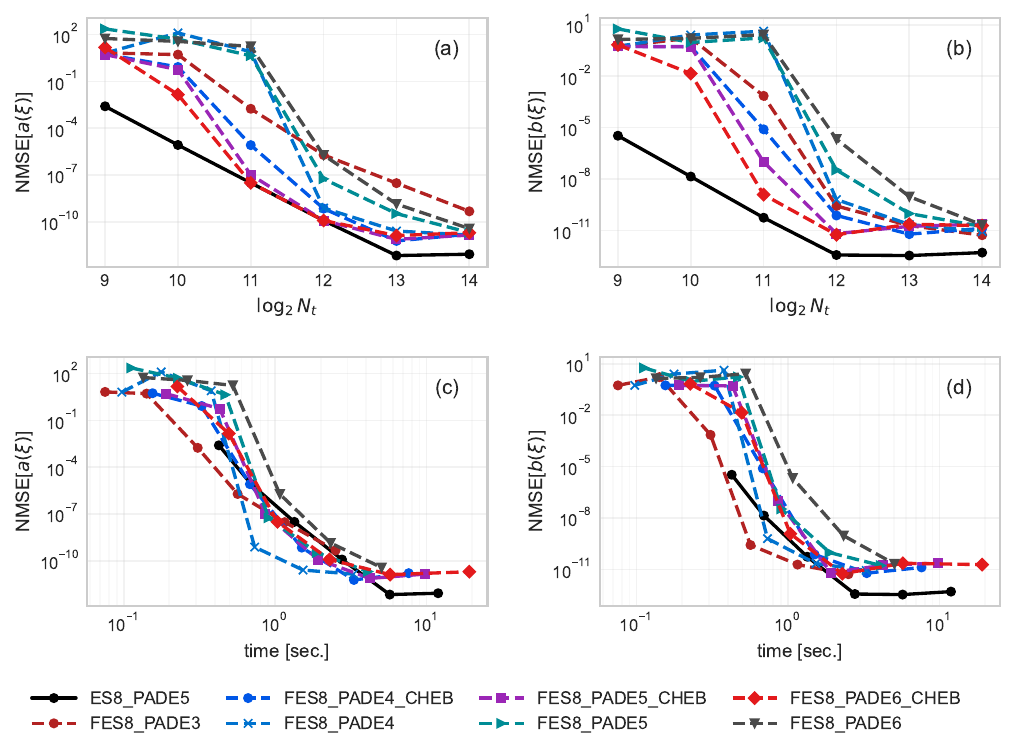}
		\caption{Continuous-spectrum errors for the chirped hyperbolic secant potential with $A=5.2$, $C=4$, and normal dispersion $\sigma=-1$. The upper panels show $\operatorname{NMSE}[a(\xi)]$ and $\operatorname{NMSE}[b(\xi)]$ as functions of the time grid size. The lower panels show the same errors as functions of the execution time. The spectral interval is $\xi\in[-20,20]$ and $N_\xi=2^{15}$.}
		\label{fig:es8_sigma_minus}
	\end{figure}
	
%	Figure~\ref{fig:es8_sigma_minus} presents the same comparison for normal dispersion, $\sigma=-1$. The behavior is similar to the anomalous-dispersion case: the Chebyshev polynomial basis improves the accuracy of the fast variants of the eighth-order scheme and makes them suitable for the computation of the continuous spectrum on the finite spectral interval. Thus, the proposed representation is not tied to a particular sign of dispersion and can be used in the continuous-spectrum setting for both cases considered here.
	
%	The numerical results demonstrate that the Chebyshev representation should be regarded as a practical way to construct fast eighth-order Pad\'e schemes. The method does not change the local eighth-order approximation itself; it changes the polynomial representation used for the fast multiplication of local transition matrices. This change is sufficient to improve the accuracy of the fast variants of the eighth-order scheme in the considered continuous-spectrum tests.

    Figures~\ref{fig:es8_sigma_plus} and~\ref{fig:es8_sigma_minus} show the results for anomalous $(\sigma=1)$ and normal $(\sigma=-1)$ dispersion. In the legends, the prefix \texttt{F} marks fast schemes and \texttt{CHEB} marks Chebyshev-based variants. The upper panels show errors versus the time-grid size and the lower panels versus execution time. The direct fast variants follow the convergence trend but are less accurate than the non-fast Pad\'e reference.
For sufficiently resolved grids, the Chebyshev-based variants substantially reduce the errors in both $a(\xi)$ and $b(\xi)$ for either sign of dispersion while retaining the fast product tree. %\textcolor{red}
{Thus, they provide a practical fast implementation of the eighth-order Pad\'e schemes on a finite spectral interval; the local eighth-order approximation remains unchanged, and only the polynomial representation used in the fast multiplication is modified.}
	
    %\clearpage
	\section{Conclusion}
	
	We proposed fast eighth-order Pad\'e schemes for the direct Zakharov--Shabat problem based on the Chebyshev polynomial representation of the spectral dependence. The local transition matrix is obtained from an eighth-order exponential scheme constructed using the Magnus expansion and is converted into rational form by diagonal Pad\'e approximation. The main difference from the direct fast Pad\'e implementation is that the spectral parameter is first restricted to a finite real interval and then represented through the Joukowski mapping and the Chebyshev polynomial basis.
	
	Numerical experiments for chirped hyperbolic secant potentials with $C=4$ and both signs of dispersion show that the direct fast variants of the eighth-order schemes may lose accuracy in the computation of the scattering coefficients. The Chebyshev-based representation improves the accuracy of these fast variants while preserving the fast product-tree construction of the global transition matrix. Therefore, the proposed approach provides a practical fast implementation of eighth-order Pad\'e schemes for computing the continuous nonlinear spectrum on a finite spectral interval.
	
	Further work will be devoted to a more detailed optimization of polynomial multiplication in the Chebyshev basis and to the extension of the approach to the computation of discrete spectral data.
	
	\section{Appendix}
	
	Here we present the expressions for the matrices $\hat{Z}_1$, $\hat{Z}_3$, $\hat{Z}_5$, $\hat{Z}_7$ included in the transition matrix~\eqref{ES8} for the ZS system~\eqref{psit}: 
	%\begin{equation}
	%	T^{(ES8)}=e^{Z(t)},\quad Z(t)=\tau Z_1(t)+\tau^3Z_3(t)+\tau^5Z_5(t) +\tau^7 Z_7(t).
	%\end{equation}
	\begin{equation}\label{Z1eZ3e}
		\hat{Z}_1=Q,\quad \hat{Z}_3=\frac{1}{24}Q^{(2)}+\frac{1}{12}\left[Q^{(1)},Q\right],
	\end{equation}
	\begin{equation}\label{Z5e}
		\hat{Z}_5=\frac{1}{1920}Q^{(4)}+\frac{1}{480}\left[Q^{(3)},Q\right]+\frac{1}{480}\left[Q^{(1)},Q^{(2)}\right]
		+\frac{1}{720}\left[\left[Q^{(2)},Q\right],Q\right]+
	\end{equation}
	$$
	+\frac{1}{240}\left[\left[Q,Q^{(1)}\right],Q^{(1)}\right]
	+\frac{1}{720}\left[Q^3,Q^{(1)}\right]+\frac{1}{240}\left[QQ^{(1)}Q,Q\right].
	$$
	\begin{equation}\label{Z7}
		\begin{array}{l}
			\displaystyle
			\hat{Z}_7 = \frac{Q Q^{(3)} Q^{(1)}}{13440} + \frac{Q^{(1)} Q^{(2)} Q^2}{40320} - \frac{(Q^{(1)})^2 Q^3}{10080} + \frac{Q Q^{(1)} Q^{(3)}}{8064}
			+ \frac{Q^{(1)} Q^3 Q^{(1)}}{3780}\\[5mm]
			\displaystyle
			- \frac{(Q^{(1)})^2 Q Q^{(1)}}{2240}  - \frac{Q^{(2)} Q^4}{30240} + \frac{Q^{(1)} Q^5}{30240}- \frac{31 Q Q^{(2)} Q^{(1)} Q}{120960} - 
			\frac{Q 
				Q^{(1)} Q^4}{6048}  + \frac{Q^{(5)} Q}{53760} \\[5mm]
			\displaystyle
			- \frac{Q^{(4)} Q^{(1)}}{53760} - \frac{Q Q^{(5)}}{53760} + \frac{Q^{6}}{322560} - \frac{Q^{(1)} Q^{(2)} Q^{(1)}}{6720} + \frac{Q Q^{(2)} 
				Q^3}{7560}+ \frac{(Q^{(1)})^3 Q}{6720} \\[5mm]
			\displaystyle
			+ \frac{Q Q^{(1)} Q Q^{(2)}}{8064} - \frac{Q (Q^{(1)})^3}{6720} + \frac{Q^3 Q^{(2)} Q}{7560} + \frac{37 Q Q^{(2)} Q Q^{(1)}}{120960} + 
			\frac{11 
				Q^{(1)} Q^2 Q^{(2)}}{120960}\\[5mm]
			\displaystyle
			+ \frac{Q^{(3)} Q^{(2)}}{16128} + \frac{23 Q^{(2)} Q^{(1)} Q^2}{120960}
			+ \frac{Q^{(1)} Q^{(3)} Q}{13440} - \frac{Q Q^{(1)} Q^2 Q^{(1)}}{2520} + \frac{Q^{(4)} Q^2}{40320} + \frac{Q^4 Q^{(1)} Q}{6048}\\[5mm] 
			\displaystyle
			- \frac{Q^5 Q^{(1)}}{30240} - \frac{Q^{(2)} Q Q^{(1)} Q}{8064} - \frac{Q^{(3)} Q^3}{40320} + \frac{Q^{(1)} Q (Q^{(1)})^2}{2240} + 
			\frac{Q^{(1)} 
				Q^{(4)}}{53760} + \frac{Q^{(2)} (Q^{(1)})^2}{13440} \\[5mm]
			\displaystyle
			- \frac{11 Q^{(2)} Q^2 Q^{(1)}}{120960} - \frac{Q (Q^{(1)})^2 Q^2}{30240} - \frac{Q^{(2)} Q^{(3)}}{16128} + \frac{Q Q^{(3)} Q^2}{13440} - 
			\frac{(Q^{(2)})^2 Q}{24192} - \frac{Q (Q^{(2)})^2}{24192}\\[5mm]
			\displaystyle
			+ \frac{Q Q^{(1)} Q Q^{(1)} Q}{7560} + \frac{Q^{(2)} Q Q^{(2)}}{12096} - \frac{37 Q^{(1)} Q Q^{(2)} Q}{120960} + \frac{(Q^{(1)})^2 
				Q^{(2)}}{13440} + \frac{31 Q Q^{(1)} Q^{(2)} Q}{120960}\\[5mm]
			\displaystyle
			- \frac{Q^{(1)} Q Q^{(3)}}{5040} + \frac{Q^{(1)} Q Q^{(1)} Q^2}{3024} - \frac{Q Q^{(4)} Q}{20160} + \frac{Q^{(3)} Q^{(1)} Q}{8064} - 
			\frac{Q^{(1)} 
				Q^2 Q^{(1)} Q}{2520} \\[5mm]
			\displaystyle
			- \frac{Q^{(3)} Q Q^{(1)}}{5040} - \frac{Q^2 Q^{(2)} Q^2}{5040}
			- \frac{Q^2 Q^{(3)} Q}{13440} - \frac{Q^2 (Q^{(1)})^2 Q}{30240} - \frac{23 Q^2 Q^{(1)} Q^{(2)}}{120960} - \frac{Q^2 Q^{(2)} 
				Q^{(1)}}{40320}\\[5mm]
			\displaystyle
			+ \frac{Q^2 Q^{(1)} Q Q^{(1)}}{3024} - \frac{Q^4 Q^{(2)}}{30240} + \frac{Q^3 Q^{(3)}}{40320} - \frac{Q^3 (Q^{(1)})^2}{10080}
			+ \frac{Q^2 Q^{(4)}}{40320} + \frac{Q^2 Q^{(1)} Q^3}{3024} - \frac{Q^3 Q^{(1)} Q^2}{3024}
		\end{array}
	\end{equation}
	
	To obtain a consistent eighth-order finite-difference scheme, it is necessary to use at least a seven-point stencil. The coefficients of the 
	finite-difference approximation of derivatives from the first order $Q^{(1)}$ to the sixth-order $Q^{(6)}$ on a seven-point uniform stencil are 
	presented in Table~\ref{tabl_deriv}.
	\begin{table}[!ht]
		\centering\small
		\caption{The coefficients of the central differences}\vspace*{2mm}
		\begin{tabular}{|c|c|c|c|c|c|c|c|}
			\hline
			Derivative & $-3$   & $-2$    & $-1$   & $0$      & $+1$    & $+2$    & $+3$\\\hline
			$1$ & $-1/60$ & $3/20$   & $-3/4$  & $0$      & $3/4$    & $-3/20$  & $1/60$\\\hline
			$2$ & $1/90$  & $-3/ 20$ & $3/2$   & $-49/18$ & $3/2$    & $-3/20$  & $1/90$\\\hline
			$3$ & $1/8$ & $-1$     & $13/8$ & $0$      & ${-13/8}$ & $1$      & $-1/8$\\\hline
			$4$ & $-1/ 6$ & $2$      & $-13/2$  & $28/3$   & $-13/2$   & $2$      & $-1/ 6$ \\\hline
			$5$ & $-1/2$  & $ 2$     & $-5/2$    & $0$      & $5/2$      & $-2$     & $1/2$ \\\hline
			$6$ & $1$     & $-6$     & $15$    & $-20$    & $15$     & $-6$     & $1$\\\hline
		\end{tabular}
		\label{tabl_deriv}
	\end{table}
	
	To compute the matrix exponential, it is convenient to use the formula  
	\begin{equation}\label{e^Zfinal}
		e^Z=c(\lambda)\sigma_0+\frac{s(\lambda)}{\lambda}Z,\quad c(\lambda)=\cosh(\lambda),\quad s(\lambda)=\sinh(\lambda),\quad  
		\lambda_\pm=\pm\sqrt{z_1^2+z_2^2+z_3^2},
	\end{equation}
	where $z_1$, $z_2$, $z_3$ are the coefficients of the expansion of the traceless matrix $Z$  ($Z_{11} = -Z_{22}$) in terms of the Pauli matrices  
	\begin{equation}\label{paulimatr}
		\sigma_0=\begin{bmatrix}1&0\\0&1\end{bmatrix}\equiv I,\quad \sigma_1=\begin{bmatrix}0&1\\1&0\end{bmatrix},\quad
		\sigma_2=\begin{bmatrix}0&-i\\i&0\end{bmatrix},\quad \sigma_3=\begin{bmatrix}1&0\\0&-1\end{bmatrix},
	\end{equation}
	namely,
	\begin{equation}\label{Zpauli}
		Z=\begin{bmatrix}Z_{11}&Z_{12}\\Z_{21}&-Z_{11}\end{bmatrix}=\begin{bmatrix}z_3&z_1-iz_2\\z_1+iz_2&-z_3\end{bmatrix}=z_1\sigma_1+z_2\sigma_2+z_3\sigma_3.
	\end{equation}
	From this, one can find the coefficients $z_1$, $z_2$, $z_3$ of the expansion in the Pauli matrices:
	\begin{equation}\label{coefZpauli}
		z_1=\frac{Z_{12} + Z_{21}}{2},\quad z_2=i\frac{Z_{12} - Z_{21}}{2},\quad z_3=\frac{Z_{11} - Z_{22}}{2}.
	\end{equation}
	
When using the diagonal Pad\'e approximation of degree $k$ for the exponential scheme%\textcolor{red}
{~\eqref{ES8}}, the transition matrix
is written as a rational function
	\begin{equation}\label{PadeGeneral}
		E_k(z)=\frac{P_k(z)}{P_k(-z)},\quad e^z=E_k(z)+O(z^{2k+1}).
	\end{equation}
%For the eighth-order scheme, one needs to use the diagonal Pad\'e approximation starting from the fourth degree. 
In the numerical experiments, approximations with $k=3,4,5,6$ are used. They have the form
    \begin{equation}\label{Pade3}
	E_3(z)=\frac{1+\frac{1}{2}z+\frac{1}{10}z^2+\frac{1}{120}z^3}{1-\frac{1}{2}z+\frac{1}{10}z^2-\frac{1}{120}z^3},
\end{equation}
\begin{equation}\label{Pade4}
	E_4(z)=\frac{1+\frac{1}{2}z+\frac{3}{28}z^2+\frac{1}{84}z^3+\frac{1}{1680}z^4}
	{1-\frac{1}{2}z+\frac{3}{28}z^2-\frac{1}{84}z^3+\frac{1}{1680}z^4}.
\end{equation}
	\begin{equation}\label{E5Pade}
		E_5(z)=
		\frac{
			1+\frac{1}{2}z+\frac{1}{9}z^2+\frac{1}{72}z^3+\frac{1}{1008}z^4+\frac{1}{30240}z^5
		}{
			1-\frac{1}{2}z+\frac{1}{9}z^2-\frac{1}{72}z^3+\frac{1}{1008}z^4-\frac{1}{30240}z^5
		},
	\end{equation}
	\begin{equation}\label{E6Pade}
		E_6(z)=
		\frac{
			1+\frac{1}{2}z+\frac{5}{44}z^2+\frac{1}{66}z^3+\frac{1}{792}z^4
			+\frac{1}{15840}z^5+\frac{1}{665280}z^6
		}{
			1-\frac{1}{2}z+\frac{5}{44}z^2-\frac{1}{66}z^3+\frac{1}{792}z^4
			-\frac{1}{15840}z^5+\frac{1}{665280}z^6
		},
	\end{equation}
%	\begin{equation}\label{E7Pade}
%		E_7(z)=
%		\frac{
%			1+\frac{1}{2}z+\frac{3}{26}z^2+\frac{5}{312}z^3+\frac{5}{3432}z^4
%			+\frac{1}{11440}z^5+\frac{1}{308880}z^6+\frac{1}{17297280}z^7
%		}{
%			1-\frac{1}{2}z+\frac{3}{26}z^2-\frac{5}{312}z^3+\frac{5}{3432}z^4
%			-\frac{1}{11440}z^5+\frac{1}{308880}z^6-\frac{1}{17297280}z^7
%		}.
%	\end{equation}
Finally, the transition matrix for the scheme based on the Pad\'e approximation can be written in a form analogous to%\textcolor{red}
{~\eqref{e^Zfinal}}
\begin{equation}\label{TPade}
		T_n(\zeta)
		=
		E_k(Z_n(\zeta))
		=c_k(x)\sigma_0+\frac{s_k(x)}{x}Z,\quad x =  i\tau\lambda.
\end{equation}

The expressions for the elements of the matrix $Z$ in~\eqref{Z5e} and, especially, in~\eqref{Z7} for the ES8 scheme~\eqref{ES8} are rather cumbersome; however, they have the form of polynomials in the spectral parameter $\zeta$. To compute the transition matrix in the form~\eqref{e^Zfinal} or~\eqref{TPade}, it is convenient to represent the Pauli expansion coefficients
$z_1$, $z_2$, $z_3$ defined by~\eqref{coefZpauli} also as polynomials in the variable $z = \tau\zeta$. Since the coefficients of these polynomials do not depend on the spectral parameter, they can be computed in advance for each value of $t$ and then used in the process of solving the ZS system for each value of the spectral parameter.

The coefficients $z_1$, $z_2$, $z_3$ of the expansion in the Pauli matrices in the case of the eighth-order scheme take the form of fifth-order polynomials
\begin{equation}\label{PauliPol8}
		z_k = c_0^k + c_1^kz + c_2^kz^2 + c_3^kz^3  + c_4^kz^4  + c_5^kz^5,\quad k = 1,2,3.
\end{equation}
Below, we give the expressions for the coefficients of these polynomials. For $z_1$ they take the form
	\begin{equation}
		\begin{array}{ll}
			c_0^1 = & \Bigl[504 q^{(4)} + 3 q^{(6)} + 967680 r - 8064 (q^{(1)})^2 r + 80 (q^{(2)})^2 r - 384 q^{(1)}  q^{(3)}r\\[5mm]
			&     - 48 q^{(4)}  r^2	+ 128 (q^{(1)})^2 r^3 + 8064 q^{(1)}  r r^{(1)} + 144 q^{(3)}  r r^{(1)} + 40320 r^{(2)}\\[5mm]
			&     - 144 (q^{(1)})^2 r^{(2)} + 144 q^{(1)}  r^{(1)} r^{(2)} + 128 q^3 ((r^{(1)})^2 + 2 r r^{(2)}) + 240 q^{(1)}  r r^{(3)}\\[5mm]
			&- 16 q^{(2)}  \bigl[-2520 + 168 r^2 - 9 q^{(1)}  r^{(1)} + 9 (r^{(1)})^2 + 5 r r^{(2)}\bigr]\\[5mm]
			& - 16 q^2 \bigl[16 q^{(2)}  r^2 + 64 q^{(1)}  r r^{(1)} - 56 r (r^{(1)})^2 + 168 r^{(2)} + 16 r^2 r^{(2)} + 3 r^{(4)}\bigr]\\[5mm]
			& + 16 q  \bigl[60480 + 3 q^{(4)}  r + 56 (q^{(1)})^2 r^2 + 504 q^{(1)}  r^{(1)} + 15 q^{(3)}  r^{(1)}\\[5mm]
			&- 64 q^{(1)}  r^2 r^{(1)} - 504 (r^{(1)})^2 + q^{(2)}  (168 r + 16 r^3 - 5 r^{(2)}) + 168 r r^{(2)}\\[5mm]
			&  + 5 (r^{(2)})^2 + 9 q^{(1)}	r^{(3)} - 24 r^{(1)} r^{(3)} + 3 r r^{(4)}\bigr]  + 504 r^{(4)}  + 3 r^{(6)}  \Bigr] / 1935360;
		\end{array}
	\end{equation}
	\begin{equation}
		\begin{array}{ll}
			c_1^1 = &i \Bigl[9 q^{(5)} - 48 q^{(3)}  (-21 + qr) - 40320 r^{(1)} + 288 (q^{(1)})^2 r^{(1)} + 64 q  q^{(2)}  r^{(1)}\\[5mm]
			& + 2688 qr r^{(1)} + 272 q^{(2)}  r r^{(1)} - 256 q^2 r^2 r^{(1)} + 208 q  r^{(1)} r^{(2)} \\[5mm]
			& + 16 q^{(1)}  \bigl[2520 - 168 qr - 13 q^{(2)}  r + 16 q^2 r^2 - 18 (r^{(1)})^2 - 17 q  r^{(2)} - 4 r r^{(2)}\bigr] \\[5mm]
			& - 1008 r^{(3)} + 48 qr r^{(3)} - 9 r^{(5)}\Bigr] / 483840;
		\end{array}
	\end{equation}
	\begin{equation}
		\begin{array}{ll}
			c_2^1 = & \Bigl[-8 q^2 r^{(2)} + 40 q  q^{(1)}  r^{(1)} + 8 q^{(2)}  (3 qr - r^2 - 21) + 24 qr r^{(2)} - 24 q  (r^{(1)})^2\\[5mm]
			& - 24 (q^{(1)})^2 r + 40 q^{(1)}  r r^{(1)} - 3 q^{(4)} - 168 r^{(2)} - 3 r^{(4)}\Bigr]  / 60480;\\[5mm]
			c_3^1 = & i \Bigl[ 3 q^{(3)} + 8 q^{(1)}  (21 - 4 qr) + 8 (-21 + 4 qr) r^{(1)} - 3 r^{(3)}\Bigr] / 30240;\\[5mm]
			c_4^1 = & -(q^{(2)} + r^{(2)}) / 3780;\quad
			c_5^1 =  i(q^{(1)} - r^{(1)}) / 1890.
		\end{array}
	\end{equation}
For $z_2$, the coefficients of the polynomial representation~\eqref{PauliPol8} take the form
	\begin{equation}
		\begin{array}{ll}
			c_0^2 = & i\Bigl[128q^3(2rr^{(2)} + (r^{(1)})^2) - 128(q^{(1)})^2r^3  - 384q^{(1)} q^{(3)} r  - 967680r\\[5mm]
			&  - 16q^2 \bigl[64q^{(1)} rr^{(1)} + 16q^{(2)} r^2 - 16r^2r^{(2)} + 56r(r^{(1)})^2 + 168r^{(2)} + 3r^{(4)}\bigr]\\[5mm]
			&  + 16q   \bigl[56(q^{(1)})^2r^2 + 64q^{(1)} r^2r^{(1)} + 504q^{(1)} r^{(1)} + 504(r^{(1)})^2 + 9q^{(1)} r^{(3)} \\[5mm]
			& + q^{(2)} (-16r^3 + 168r - 5r^{(2)}) + 15q^{(3)} r^{(1)} + 3q^{(4)} r - 168rr^{(2)} - 3rr^{(4)}\\[5mm]
			&   + 24r^{(1)}r^{(3)} - 5(r^{(2)})^2 + 60480\bigr] - 8064(q^{(1)})^2r - 144(q^{(1)})^2r^{(2)}\\[5mm]
			& + 16q^{(2)} \bigl[9(q^{(1)} r^{(1)} + (r^{(1)})^2 + 280) + 168r^2 + 5rr^{(2)}\bigr]- 8064q^{(1)} rr^{(1)} \\[5mm]
			&  - 240q^{(1)} rr^{(3)} - 144q^{(1)} r^{(1)}r^{(2)} + 80(q^{(2)})^2r - 144q^{(3)} rr^{(1)}\\[5mm] 
			&  + 48q^{(4)} r^2 + 504q^{(4)} + 3q^{(6)} - 40320r^{(2)} - 504r^{(4)} - 3r^{(6)}\Bigr]/1935360;
		\end{array}
	\end{equation}
	\begin{equation}
		\begin{array}{ll}
			c_1^2 = & \Bigl[-9 q^{(5)} + 48 q^{(3)}  (-21 + qr) - 40320 r^{(1)} - 288 (q^{(1)})^2 r^{(1)} - 64 q  q^{(2)}  r^{(1)} \\[5mm]
			&+ 2688 qr r^{(1)} + 272 q^{(2)}  r r^{(1)} - 256 q^2 r^2 r^{(1)} + 208 q  r^{(1)} r^{(2)} \\[5mm]
			&- 16 q^{(1)}  \bigl[2520 - 168 qr - 13 q^{(2)}  r + 16 q^2 r^2 + 18 (r^{(1)})^2 - 17 q  r^{(2)} + 4 r r^{(2)}\bigr] \\[5mm]
			&	- 1008 r^{(3)} + 48 qr r^{(3)} - 9 r^{(5)}\Bigr] / 483840;
		\end{array}
	\end{equation}
	\begin{equation}
		\begin{array}{ll}
			c_2^2 = & -i \Bigl[3 q^{(4)} + 24 (q^{(1)})^2 r - 8 q^{(2)}  (-21 + 3 qr + r^2) - 40 q  q^{(1)}  r^{(1)} \\[5mm]
			&+ 40 q^{(1)}  r r^{(1)} - 24 q  (r^{(1)})^2 - 168 r^{(2)} + 8 q^2 r^{(2)} + 24 qr r^{(2)} - 3 r^{(4)}\Bigr] / 60480;\\[5mm]
			c_3^2 = & \Bigl[-3 q^{(3)} + 8 q^{(1)}  (-21 + 4 qr) - 168 r^{(1)} + 32 qr r^{(1)} - 3 r^{(3)}\Bigr] / 30240;\\[5mm]
			c_4^2 = & i (r^{(2)} - q^{(2)} ) / 3780;\quad
			c_5^2 = -(q^{(1)} + r^{(1)}) / 1890.
		\end{array}
	\end{equation}
For $z_3$, the coefficients of the polynomial representation~\eqref{PauliPol8} take the form
	\begin{equation}
		\begin{array}{ll}
			c_0^3 = & \Bigl[-256 q^3 r^2 r^{(1)}  + 2688 q^2 r r^{(1)} + 48 q^2 r r^{(3)} + 208 q^2 r^{(1)} r^{(2)}  \\[5mm]
			& + q^{(1)}\bigl[256 q^2 r^3 - 16 r^2 (168 q + 13 q^{(2)} ) - 160 r (q  r^{(2)} - 252) \\[5mm]
			& + 9 [-32 q  (r^{(1)})^2 + 112 r^{(2)} + r^{(4)}]\bigr] + 160 q  q^{(2)}  r r^{(1)} \\[5mm]
			& - 6 q^{(3)}  (8 q  r^2 - 168 r - 5 r^{(2)}) - 40320 q  r^{(1)} - 1008 q  r^{(3)} - 9 q  r^{(5)}\\[5mm]
			&  + 288 (q^{(1)})^2 r r^{(1)} - 1008 q^{(2)}  r^{(1)} - 30 q^{(2)}  r^{(3)} - 9 q^{(4)}  r^{(1)} + 9 q^{(5)}  r\Bigr] / 483840;
		\end{array}
	\end{equation}
	\begin{equation}
		\begin{array}{ll}
			c_1^3 = & -i \Bigl[60480 - 3 q^{(4)}  r + 8 (q^{(1)})^2 r^2 + 1008 q^{(1)}  r^{(1)} + 24 q^{(3)}  r^{(1)} - 112 q  q^{(1)}  r r^{(1)} 
			\\[5mm]
			&+ 8 q^2 (r^{(1)})^2 + 2 q^{(2)}  (-84 r + 8 q  r^2 - 5 r^{(2)}) - 168 q  r^{(2)} + 16 q^2 r r^{(2)} \\[5mm]
			&+ 24 q^{(1)}  r^{(3)} - 3 q  r^{(4)}\Bigr] / 60480;
		\end{array}
	\end{equation}
	\begin{equation}
		\begin{array}{ll}
			c_2^3 = & \Bigl[3 q^{(3)}  r - 168 q  r^{(1)} + 11 q^{(2)}  r^{(1)} + 32 q^2 r r^{(1)} \\[5mm]
			& + q^{(1)}  (168 r - 32 q  r^2 - 11 r^{(2)}) - 3 q  r^{(3)}\Bigr] / 30240;\\[5mm]
			c_3^3 = & i \Bigl[q^{(2)}  r - 8 q^{(1)}  r^{(1)} + q  r^{(2)}\Bigr] / 3780;\quad
			c_4^3 = (q^{(1)}  r - q  r^{(1)}) / 1890;\quad
			c_5^3 = 0.
		\end{array}
	\end{equation}

	\bibliographystyle{unsrt}
	\bibliography{references}

\end{document}